\documentclass[12pt]{amsart}
\usepackage{amssymb,amscd}
\usepackage{graphicx,wrapfig,graphics}
\newcommand{\N}{{\mathbb N}}
\newcommand{\T}{{\mathbb T}}
\newcommand{\Z}{{\mathbb Z}}
\newcommand{\Q}{{\mathbb Q}}
\newcommand{\R}{{\mathbb R}}
\newcommand{\C}{{\mathbb C}}
\newcommand{\St}{\frak S\frak t\frak a\frak b\frak l\frak e}
\newcommand{\Un}{\frak U\frak n\frak s\frak t\frak a\frak b\frak l\frak e}
\newtheorem{theorem}{Theorem}
\newtheorem{lemma}[theorem]{Lemma}

\newtheorem{corollary}[theorem]{Corollary}
\newtheorem{conjecture}[theorem]{Conjecture}
\newtheorem{example}[theorem]{Example}

\newcommand{\holom}{{\mathcal H}}
\newcommand{\laurent}{{\mathcal L}}
\newcommand{\polynom}{{\mathbb Z}[x]}
\newcommand{\zsup}{\Z_{ac}(x)}
\newcommand{\zrightsup}{{\mathbb Z}_{ac}[x)}

\newcommand{\Tbiinfinite}{{\mathbb T}(x)}

\title
{Self Duality and Codings for Expansive Group Automorphisms}
\address{Faculty of Mathematics, University of North Texas, Denton, Texas, USA}
\address{Delft University, Faculty of Electrical Engineering, Mathematics and
Information Technology, P.O.Box 5031, 2600 GA Delft, Netherlands}
\author {Alex Clark}\author{Robbert Fokkink}

\begin{document}

\maketitle
\begin{abstract}
Lind and Schmidt have shown that the homoclinic group of a cyclic
$\Z^k$ algebraic dynamical system is isomorphic to the dual of the
phase group. We show that this duality result is part of an exact
sequence if $k=1$. The exact sequence is a well known algebraic
object, which has been applied by Schmidt in his work on rigidity.
We show that it can be derived from dynamical considerations only.
The constructions naturally lead to an almost $1-1$-coding of
certain Pisot automorphisms by their associated $\beta$-shift,
generalizing similar results for Pisot automorphisms of the torus.
\end{abstract}

\section{Introduction}
In this paper we study discrete systems in which the phase space
is a topological group and the action is induced by an expansive
automorphism.

\medbreak

A common practice in algebra is to embed a ring $R$ as a discrete
subring of a locally compact ring $K$ such that $K/R$ is compact
and is equal to the Pontryagin dual of $R$ as an additive group.
Standard examples of such embeddings are: $R$ is the ring of
integers of a number field $\Q(\alpha)$ and $K$ is the product of
the archimedean completions of $\Q(\alpha)$; or, $R$ is a number
field and $K$ is the ad\`ele ring. Pontryagin duality concerns
groups and not rings, so the multiplicative structure on $R$ is
lost in $K/R$. However it can be regained by considering $R$ and
$K/R$ as dynamical systems. For instance, if $R$ is equal to
$\Z[\beta,\beta^{-1}]$ then the multiplicative structure on $R$ is
determined by the action $x\mapsto\beta x$ and so $K/R$ is a
compact group with a dual action; i.e., it is an algebraic
dynamical system.

\medbreak
In this paper we build up the exact sequence
\[
\begin{CD}
0 @>{}>> R @>{}>> K @>>> G @>>> 0
\end{CD}
\]
by dynamical considerations, starting from $G$ rather than $R$:
$K$ is the product of the stable and unstable group of $G$,
endowed with suitable topologies, and $R$ is the homoclinic group.
Our paper was motivated by work of Lind and Schmidt
\cite{Lindschmidt}, who showed that the homoclinic group of an
algebraic dynamical system is isomorphic to the Pontryagin dual,
and by work of Schmidt on symbolic codings \cite{Schmidt}. We only
consider $\Z$-actions. Lind and Schmidt obtained their results in
a more general case of ergodic $\Z^k$ actions on compact abelian
groups.

\medbreak

We give two applications. We interpret a number theoretic result
of Pisot and Vijayaraghavan in dynamical terms and we find an
almost $1-1$-coding of a class of Pisot automorphisms. Such
codings were known to exist in the case the underlying group is a
torus (see \cite{Schmidt,Sidorovvershik,Sidorov1}), but our
constructions apply to a more general class of groups. Our
analysis of the coding also leads to an interpretation of the adic
transformation associated to a $\beta$-shift as the return map of
a naturally related flow.

\section{Notation}

We shall study certain rings of polynomials and power series that
are a little unusual, so our notation is a little unusual. If $R$
is a ring and $x$ is an indeterminate, then $R[x]$ denotes the
ring of Laurent polynomials $r_0x^n+r_1x^{n+1}+\ldots+ r_kx^{n+k}$
with $r_i\in R$, $k\in\N$ and $n\in\Z$. $R[x)$ denotes the ring of
power series over $R$ of the form
$r_0x^n+r_1x^{n+1}+\ldots+r_kx^{n+k}+\ldots$ for $k=0,1,2,\ldots$.
Similarly, $R(x]$ denotes the ring of power series over $R$ of the
form $\ldots+r_kx^{n-k}+\ldots+r_1x^{n-1}+r_0x^{n}$. The ring of
two-sided or Laurent power series over $R$ is denoted by $R(x)$.

\section{Cyclic algebraic dynamical systems}

A dynamical system is a group action on a topological space
$\mathcal{G}\times X\to X$. If the phase space $X$ is an abelian
topological group and if the action $x\mapsto g\cdot x$ is an
automorphism, then the system is called an \textit{algebraic
dynamical system}. In this case, $X$ is a $\Z[\mathcal{G}]$
module. These systems have been well studied and more information
on them can be found in Schmidt's excellent monograph~\cite{mono}.
In our case, $\mathcal{G}$ is $\Z$ and $X$ is a compact connected
abelian group denoted $\Gamma$ and the group ring
$\Z[\mathcal{G}]$ is isomorphic to the ring $\polynom$ of Laurent
polynomials $a_0x^n+\ldots +a_kx^{n+k}$ for $n,a_i\in\Z$ and
$k\in\N$.

\medbreak We denote the dual group of the phase space $\Gamma$ by
$G$. Then $\Gamma$ is a $\polynom$-module, as is $G$ under the
dual action. Since $\Gamma$ is connected, $G$ is torsion-free.
Hence, the annihilator of $G$ in $\polynom$ is a principal ideal,
generated by a primitive Laurent polynomial $f$. If we require
that $f=a_0+\ldots +a_n x^n$ be a polynomial; i.e., that all
powers of $x$ be positive, that $f(0)\not=0$ and that $a_n>0$,
then $f$ is uniquely determined. This polynomial $f$ is the
\textit{associated polynomial} of the dynamical system on
$\Gamma$. \medbreak

The action of $\polynom$ on $\Gamma$ is generated by the
automorphism $\gamma\mapsto x\cdot \gamma$. We denote this
automorphism by $\alpha$ and we denote the algebraic dynamical
system by $(\Gamma,\alpha)$. The stable group
$\Omega^+\subset\Gamma$ is defined as
\[
\Omega^+=\{g\in\Gamma\colon \alpha^n(g)\to 0\text{ as
}n\to\infty\}
\]
The unstable group $\Omega^-$ is defined likewise for
$n\to-\infty$. The homoclinic group is then
$\Omega=\Omega^+\cap\Omega^-$. An algebraic dynamical system
$(\Gamma,\alpha)$ is \textit{expansive} if there exists a
neighborhood $U$ of the identity $e\in G$ such that
\[
\bigcap\left\{\alpha^n(U)\colon n\in\Z\right\}=\{e\}
\]
\medbreak A polynomial is \textit{hyperbolic} if all its roots in
$\C$ are off the unit circle $|z|=1$. The following classification
theorem of expansive systems has apparently been proved first by
Aoki and Dateyama~\cite{Aoki}.

\begin{theorem}\label{aokidateyama}

An algebraic system $(\Gamma,\alpha)$ is expansive if and only if
the associated polynomial is hyperbolic and the dual group is a
Noetherian $\polynom$-module.
\end{theorem}

An algebraic system on $\Gamma$ is \textit{cyclic} or
\textit{principal} if the Pontryagin dual $G$ is a cyclic
$\polynom$-module. We shall only consider algebraic systems that
are expansive and cyclic. In this case $G$ is equivalent to the
quotient module $\polynom/(f)$ as described by the exact sequence.
\begin{equation}
\begin{CD}
\label{exact1} 0 @>{}>> \polynom @>{g\to f\cdot g}>> \polynom@>>>
G @>>> 0
\end{CD}
\end{equation}
Let $\T$ denote the circle group. The Pontryagin dual of
$\polynom$ is the group $\Tbiinfinite$ of formal power series
$\sum_{n\in\Z}t_nx^n$ with coefficients $t_n\in\T$. The Pontryagin
dual of sequence \ref{exact1} is
\begin{equation}
\begin{CD}
\label{exact2} 0 @<{}<< \Tbiinfinite @<{f\cdot g\leftarrow g}<<
\Tbiinfinite@<<< \Gamma @<<< 0
\end{CD}
\end{equation}
Some caution is required. If the action of the indeterminate $x$
on $G$ is given by multiplication by $x$, then the adjoint action
on $\Tbiinfinite$ is given by multiplication by $x^{-1}$. We have
to choose for which of the dual groups $G$ and $\Gamma$ the action
is by multiplication of $x$. We decide that the action of
$\polynom$ on $\Gamma$ is by multiplication of $x$ and the action
on $G$ is by multiplication of $x^{-1}$.

\section{Almost convergent series}

The interest in the homoclinic group arose from ideas of Vershik
on symbolic codings. Homoclinic points have exponential decay;
i.e., $\lim_{|n|\to\infty}f^{n}(x)=0$ exponentially, so that if
$(a_n)$ is a bounded sequence of integers, then $\sum_{n\in\Z}a_n
f^n(x)$ converges. This defines a conjugation between the shift on
the symbolic sequences and the action of $f$ on the phase group.
In Vershik's approach, the symbolic sequences are the integral
sequences in $l^{\infty}$, the Banach space of bounded sequences.
It turns out that the homoclinic group corresponds to the subset
of sequences for which $a_n=0$ for all but finitely many $n$. The
stable group and the unstable group correspond to sequences that
have a tail of zeroes. Lind and Schmidt analyze the homoclinic
group. For our purposes, showing that the stable and unstable
group are rings, this setting is a little inconvenient. The
multiplication of two stable elements $(a_n),(b_n)$ is well
defined: take the product of the generating power series $\sum a_n
x^n$ and $\sum b_n x^n$. However, the coefficients of this product
are not necessarily bounded. Therefore we replace $l^\infty$ by a
slightly larger space that is closed under the multiplication of
stable elements, and we call this larger space the space of almost
convergent power series. \medbreak The outer radius $R$ of a
Laurent series $\sum_{-\infty}^\infty c_n x^n$ is defined as
$R=\liminf_{n\to\infty} \sqrt[n]{|c_n|}$. The inner radius $r$ is
defined as $r=\limsup_{n\to -\infty} 1/{\sqrt[n]{|c_{n}|}}$. Let
$\holom$ be the set of Laurent series with complex coefficients
such that $r<1<R$. Then $\holom$ represents the space of analytic
functions that are defined on some domain around the unit circle.
Thus $\holom$ is a ring. Let $\laurent$ be the set of Laurent
series with complex coefficients such that $r\leq 1\leq R$. We say
that $\laurent$ is the space of \textit{almost convergent series}.
By identifying $\sum c_ix^i$ with the sequence $(c_i)\in\C^\Z$, we
endow $\laurent$ with the product topology, and we call this the
\textit{weak topology}. \medbreak

It is obvious that $\laurent$ is a $\polynom$-module. In
accordance with the action on $\Tbiinfinite$, we define the action
of $x$ on $\laurent$ as multiplication by $x$. It turns out that
this action extends to $\holom$.

\begin{lemma}\label{module}
$\laurent$ is an $\holom$-module.
\end{lemma}

\begin{proof}
Suppose that $g=\sum_{-\infty}^{\infty} c_n x^n$ is an almost
convergent series. Then $g^+=\sum_{n> 0} c_n x^n$ converges on
$|x|<1$ and $g^-=\sum_{n\leq 0} c_n x^n$ converges on $|x|>1$.
Since an element $h\in \holom$ has no poles on the unit circle,
$h\cdot g^+$ is analytic on a region $r<|x|<1$, while $h\cdot g^-$
is analytic on a region $1<|x|<R$. Thus, the outer radius of
$f\cdot(g^-+g^+)$ is $\geq 1$ and the lower radius is $\leq 1$. In
particular, $h\cdot g\in\laurent$.
\end{proof}

We shall denote almost convergent power series by the subscript
$ac$. For example, $\R_{ac}(x)$ denotes the ring of almost
convergent Laurent series with real coefficients. We shall denote
convergent power series for which $r<1<R$ by the subscript $c$.
For example, $\R_c(x)$ denotes the ring of Laurent series with
real coefficients, convergent on an annulus around the unit
circle. By Lemma~\ref{module}, $\R_{ac}(x)$ is a $\R_c(x)$-module.
Note that hyperbolic polynomials are units in $\R_c(x)$.

\begin{lemma}\label{intersect}
Suppose that $f\in\polynom$ is hyperbolic. Then $f\cdot\zsup\cap
\zrightsup=f\cdot\zrightsup$.
\end{lemma}
\begin{proof}
Suppose that $m^+\in f\cdot\zsup\cap \zrightsup$, so that
$m^+=f\cdot m$ for some $m\in\zsup$. Since $m^+$ is holomorphic on
the punctured disc and since $\frac 1 f$ is meromorphic without a
pole on the unit circle, $m=\frac 1 f \cdot m^+$ has inner radius
$r<1$. Since $m\in\zsup$, this implies that $m=\sum a_i x^i$ and
that $a_i=0$ for sufficiently small index.
\end{proof}

The projection $\R_{ac}(x)\to\T(x)$ is defined by reducing the
coefficients modulo~$1$. It obviously is a surjection. We have the
following commutative diagram of continuous group homomorphisms

\begin{equation}
\begin{CD}
0 @<<{}< \R_{ac}(x) @<<{f\cdot g\leftarrow g}<
\R_{ac}(x)@<<< 0 \\ @.@VVV@VVV@.\\
0 @<<{}< \Tbiinfinite @<<{f\cdot g\leftarrow g}< \Tbiinfinite@<<<
\Gamma@<<<0 \label{exact3}
\end{CD}
\end{equation}
in which the top row is invertible since $f$ is a unit in
$\R_c(x)$. The image of $\Gamma$ in $\T(x)$ is equal to the
projection of $\frac 1 f\cdot \Z_{ac}(x)$ on $\T(x)$. More
specifically
\begin{equation}
\label{gamma} \Gamma\cong\frac 1
f\cdot\Z_{ac}(x)/\Z_{ac}(x)\cong\Z_{ac}(x)/f\cdot\Z_{ac}(x)
\end{equation}
as algebraic $\polynom$ modules. We identify $\Gamma$ with
$\Z_{ac}(x)/f\cdot\Z_{ac}(x)$ and denote this quotient by
$\Z_{ac}(x)/(f)$. \medbreak

\begin{lemma}
The stable group $\Omega^+\subset\Gamma$ is algebraically
isomorphic to $\Z_{ac}[x)/(f)$.
\end{lemma}
\begin{proof}
A power series $g=\sum r_ix^i$ in $\frac 1f \cdot \Z_{ac}(x)$
projects onto a stable element of $\T(x)$ if and only if
$r_i\text{ mod }1$ converges to $0$ as $i\to -\infty$. Without
changing the projection of $g$ onto $\T(x)$ we may choose the
$r_i$ such that $\lim_{i\to-\infty}r_i=0$. If $f\cdot g=\sum
s_ix^i$ then, since $f$ is a polynomial, $\lim_{i\to-\infty}s_i=0$
as well. Since $f\cdot g\in\Z_{ac}(x)$ the coefficients $s_i$ are
integers and so $s_i=0$ for sufficiently small $i$. In other
words, $f\cdot g\in\Z_{ac}[x)$ and by Lemma~\ref{intersect} we see
that
\begin{equation}
\Omega^+\cong \frac 1f \cdot \Z_{ac}[x)/\Z_{ac}[x)\cong
\Z_{ac}[x)/f\cdot \Z_{ac}[x).
\end{equation}
\end{proof}

The following corollary is a special case of a result of Lind and
Schmidt~\cite{Lindschmidt} characterizing the homoclinic
group, which they have derived for certain ergodic $\Z^k$-actions.

\begin{corollary}\label{homoclinic=dual}
The homoclinic group is isomorphic to the dual group. Even more
so, the two are equivalent as algebraic dynamical systems.
\end{corollary}
\begin{proof}
$\Omega=\Omega^-\cap\Omega^+=(\Z_{ac}(x]\cap\Z_{ac}[x))/(f)=\Z_{ac}[x]/(f)=\polynom/(f).$
Not only is the homoclinic group isomorphic to the dual group, but
it is equivalent as a dynamical system as well. The action on the
dual group is multiplication by $x^{-1}$. The action on the
homoclinic group is multiplication by $x$.
\end{proof}
Notice that this representation of the homoclinic group also shows
that it possesses a \textit{fundamental} homoclinic point; that
is, a homoclinic point whose iterates span the homoclinic group.
So we have found good algebraic descriptions of $\Gamma$ and its
dynamically defined subgroups. In particular, we have an exact
sequence of algebraic groups
\begin{equation}
\begin{CD}
0 @>>> \Omega @>{x\to(-x,x)}>> \Omega^-\times\Omega^+ @>{(x,y)\to
x+y}>> \Gamma@>>> 0
\end{CD}
\end{equation}
To see that this sequence is exact, observe that if $x+y=0$ for
$(x,y)\in\Omega^-\times\Omega^+$ then $(x,y)$ can be represented
by $(g,h)\in\Z_{ac}(x]\times\Z_{ac}[x)$ such that $g+h=f\cdot k$
for some $k\in\Z_{ac}(x)$. Let $k_1$ be the sum of the
non-positive powers of $k$ and let $k_2$ be sum of the positive
powers. Then $g-f\cdot k_1=-(h-f\cdot k_2)$ and $g-f\cdot
k_1\in\Z_{ac}(x]$ while $g-f\cdot k_2\in\Z_{ac}[x)$, So
$p=g-f\cdot k_1$ is a polynomial and $(g,h)=(p,-p)\text{ mod
}(f)$. \bigbreak Now we have to find good topological descriptions
as well. $\Gamma$ is not topologically isomorphic to
$\Z_{ac}(x)/(f)$ if this quotient is endowed with the quotient
topology of $\Z_{ac}(x)$. The problem is that $g\mapsto f\cdot g$
is algebraically invertible, but the inverse is not continuous in
the weak topology. To remedy this, we shall endow $\Z_{ac}(x)$
with a topology that is stronger than the weak topology.

\section{Strong topologies on the stable group and the unstable group}

For a natural number $N$ define the subset $B(N)\subset\Z_{ac}(x)$
by
\[
B(N)=\left\{\ \sum a_i x^i\in\zsup\colon 0\leq a_i\leq N\text{ for
all }i\right\},
\]

\begin{lemma}\label{symbolic}
For every $f\in\polynom$ there exists an $N$ such that
$B(N)+f\cdot\zsup=\zsup$.
\end{lemma}
\begin{proof}
Let $N$ be equal to the sum of the absolute values of the
coefficients of $f$. Suppose that $g\in\zsup$ and that $\frac 1
f\cdot g=\sum r_ix^i$ with $r_i\in\R$. For $h=\sum\lfloor
r_i\rfloor x^i$ one verifies that $g-f\cdot h\in B(N-1)$
and that $h\in\Z_{ac}(x)$.
\end{proof}

By this lemma, the projection $B(N)\to\Z_{ac}(x)/(f)$ is onto
provided that $N$ is sufficiently large. We define the
\textit{strong topology} on $\Z_{ac}(x)/(f)$ as the quotient
topology induced by $B(N)$ for some $N$ that is sufficiently
large. In particular, $\Z_{ac}(x)/(f)$ is compact. Moreover, the
choice of $N$ does not alter this topology.

\begin{theorem}\label{same}
$\Gamma$ is isomorphic to $\Z_{ac}(x)/(f)$ with the strong
topology as a $\polynom$-module.
\end{theorem}
\begin{proof}
It suffices to show that the inverse map $g\mapsto \frac 1 {\bar
f} \cdot g$ restricted to $B(N)$ is continuous. Since $f$ is a
hyperbolic polynomial, it has an inverse $\frac 1 {\bar f} =\sum
q_i x^i$ in $\R_{c}(x)$. For every $\epsilon>0$ there exists an
$n$ such that $\sum_{|i|>n}|q_i|<\epsilon/N$. Let $0\in U\subset
B(N)$ be the neighborhood defined by $U=\left\{\sum a_ix^i\colon
a_i=0 \text{ if } |i|\leq n\right\}$. Then all elements $\sum r_i
x^i\in\frac 1{\bar f}\cdot U$ have coefficient~$|r_0|<\epsilon$.
By choosing a larger $n$, we can restrict arbitrarily many
coefficients $r_k$ by $\epsilon$.
\end{proof}

The stable group $\Omega^+\subset \Gamma$ is a dense subgroup,
provided $f$ is not a unit in $\Z_{ac}[x)$. So, in general, the
stable group is not locally compact. To remedy that, we endow the
stable group with a stronger topology. For a natural number $N$
define the subset $B_k(N)\subset\Z_{ac}[x)$ by
\[
B_k(N)=\left\{\ \sum a_i x^i\colon  |a_i|\leq N\text{ for all
}i\text{ and } a_i=0\text{ for }i<-k\right\},
\]
Clearly, $B_k(N)$ is a compact set in the weak topology. Let
$B_\infty(N)$ be the union of all $B_k(N)$.

\begin{lemma}\label{symbolic2}
For every $f\in\polynom$ there exists an $N$ such that
$B_\infty(N)+f\cdot\Z_{ac}[x)=\Z_{ac}[x)$.
\end{lemma}
\begin{proof}
Let $N$ be the sum of all the absolute values of the coefficients
of $f$. Suppose that $g\in\Z_{ac}[x)$ and that $\frac 1 f\cdot
g=\sum r_ix^i$ with $r_i\in\R$. Then $\lim_{i\to-\infty}|r_i|=0$.
Let $[r_i]$ be the integer that is nearest to $r_i$. For $h=\sum
[r_i] x^i$ one verifies that $g-f\cdot h\in B_\infty(N)$ and that
$h\in\Z_{ac}[x)$.
\end{proof}

Let $N$ be a sufficiently large integer so that $B_\infty(N)$
projects onto $\Omega^+$. We define the strong topology on
$\Omega^+$ by the filtration $B_k(N)$. In particular,
$U\subset\Omega^+$ is closed if and only if its preimage in
$\tilde U\subset B_\infty(N)$ intersects each $B_k(N)$ in a closed
subset. By symmetry, we endow $\Omega^-$ with the strong topology,
and we endow $\Omega$ with the discrete topology.

\begin{theorem}[Unfolding]
The sequence of locally compact abelian groups
\begin{equation}
\label{unfolding}
\begin{CD}
0 @>>> \Omega @>{x\to(-x,x)}>> \Omega^-\times\Omega^+ @>{(x,y)\to
x+y}>> \Gamma@>>> 0
\end{CD}
\end{equation}
 is a self-dual exact sequence.
\end{theorem}
\begin{proof}
We have already seen that the sequence is algebraically exact. It
is easy to verify that all maps are continuous. To prove that the
maps are open, it suffices to show that $\Omega$ embeds as a
discrete subgroup of $\Omega^-\times\Omega^+$. Since the image of
$\Omega$ is equal to the kernel of the continuous projection
$\Omega^-\times\Omega^+\to \Gamma$, this image is closed, hence
locally compact. Since $\Omega$ is countable, this image is
countable. A countable and locally compact group is necessarily
discrete.

By Corollary~\ref{homoclinic=dual} the homoclinic group $\Omega$
is isomorphic to the dual of $\Gamma$. So to establish self
duality of the sequence, we need to show that $\Omega^+$ and
$\Omega^-$ are self-dual. We claim that the dual of $\Omega^+$ is
stable under the adjoint action. To see this, note that all
elements of $\Omega^+$ converge to $0$ under iteration of $x$. By
continuity, a character of $\Omega^+$ converges to the zero
character under the adjoint action. Hence the dual group of
$\Omega^+$ is stable under the adjoint action and, by symmetry,
the dual group of $\Omega^-$ is unstable.

Since $\Omega$ embeds discretely, both in $\Omega^-\times\Omega^+$
and in its dual, $\Omega^-\times\Omega^+$ is locally isomorphic to
its dual. Let $U\subset \Omega^-\times \Omega^+$ be a neighborhood
for which the local isomorphism is defined and let $V$ be its
image. The stable elements of $V$ are in the dual of $\Omega^+$.
We extend the local isomorphism to an isomorphism between
$\Omega^+$ and its dual, as follows. For a stable element
$s\in\Omega^+$, let $x^n\cdot s $ be the first element in $U$
under forward iteration of $x$ and let $t$ be its image in $V$.
Define the image of $s$ as $x^{-n}\cdot t$, where $x$ now denotes
the adjoint action on the dual group. This defines a homomorphism
of $\Omega^+$ to its dual group and its inverse can be defined in
the same way. So $\Omega^+$ is self-dual and by symmetry so is
$\Omega^-$.
\end{proof}

In the next section we describe $\Omega^+$ and $\Omega^-$
by algebraic methods.

\section{Stable prime divisors}

We do not use our special notation for rings of power series when
we add determinates: $\Z[\beta]$ denotes the ring generated by
$\Z$ and $\beta$, as usual. It is different from
$\Z[\beta,\beta^{-1}]$. We stick to our notation for
indeterminates, and so our $\polynom$ is usually denoted by
$\Z[x,x^{-1}]$. \medbreak In this section, we first assume that
the associated polynomial $f\in\polynom$ is \textit{irreducible}.
The case of a reducible $f$ follows from the Chinese remainder
theorem, which we deal with it at the end of this section. So for
now, the homoclinic group $\polynom/(f)$ is a ring without zero
divisors. Let $Q$ be its field of fractions and let $K_P$ be the
completion of $Q$ with respect to a prime divisor $P$. The
standard method to discretize $\polynom/(f)$ is by embedding it in
a product of $K_P$ for a certain choice of prime divisors. For
instance, if $f$ is a monic polynomial with unit constant
coefficient, then all elements of $\polynom/(f)$ are algebraic
integers and $P$ ranges over all archimedean prime divisors. In
the general case, one just needs a few prime-divisors more.
\medbreak From now on, let $\beta$ denote a root of $f$, so that
$\polynom/(f)$ is isomorphic to $\Z[\beta,\beta^{-1}]$ and the
field of fractions is isomorphic to $\Q(\beta)$. We say that a
prime-divisor $P$ of a number field $\Q(\beta)$ is \textit{stable}
if and only if $\lim_{n\to\infty} \beta^n=0$ in the $P$-topology.
Likewise, when $\lim_{n\to - \infty}\beta^{n}=0$, we say $P$ is
\textit{unstable}. Since the associated polynomial $f$ is
hyperbolic, all archimedean prime-divisors are either stable or
unstable. If $P$ is a stable prime divisor, then we say that
$\beta$ is a stable root with respect to $P$. More specifically,
we say that $\beta$ is $P$-stable if and only if $\beta$ is an
element of $O_P\subset \Q(\beta)$, the ring of integers with
respect to $P$. \medbreak We denote the set of all stable
prime-divisors of $\Q(\beta)$, archimedean or non-archimedean, by
$\St$. It is a finite set. By the Newton polygon, the
non-archimedean prime-divisors in $\St$ are extensions of the
rational primes that divide $f(0)$, the constant coefficient.

\begin{lemma}\label{unstableinverse}
Suppose that $a^{-1}\in O_P$. Then
$1/(a-x)\in O_P[x)$.
\end{lemma}
\begin{proof}
The inverse of $a-x$ is $\sum_{n=0}^\infty a^{-1-n} x^n$.
\end{proof}

Let $K_P$ be the completion of $\Q(\beta)$ with respect to $P$. If
$P$ is stable, then the evaluation map $\sum a_nx^n\mapsto \sum
a_n \beta^n$ induces a homomorphism $e_P\colon O_P[x)\to K_P$.
Since $\zrightsup\subset O_P[x)$ whenever $P$ is stable, the
product of the evaluation maps defines a homomorphism from
$\zrightsup$ to the product of $K_P$, where $P$ ranges over the
stable non-archimedean prime-divisors. Since the radius of
convergence of a power series in $\zrightsup$ is at least $1$, the
evaluation map remains well defined for stable archimedean prime
divisors as well. Moreover, it is continuous with respect to the
strong topology. We are going to show that the product of the
evaluation maps $e_P$ for $P\in\St$ gives a surjective
homomorphism with kernel $(f)$.

\begin{lemma}\label{stableinverse}
Let $P$ be a non-archimedean stable prime-divisor. If $e_{P}(g)=0$
for $g\in O_P[x)$, then $g$ is divisible by $\beta - x$.
\end{lemma}
\begin{proof}
Let $|\cdot|\in P$ be a valuation and let $g=\sum a_nx^n$ be a
power series such that $e_P(g)=0$. Then the partial sums
$s_N=\sum_{n<N} a_n \beta^n$ converge to zero in the $P$-topology
and have additive inverse $\sum_{n\geq N} a_n \beta^n$. Since
$|a_n|\leq 1$ and since $P$ is non-archimedean,
$|s_m|=|\sum_{n\geq m} a_n \beta^n|\leq |\beta|^m$. If $h=\sum
(s_n/{\beta^n}) x^{n-1}$, then $|s_n/{\beta^n}|\leq 1$ and
$g=(\beta -x)h$.
\end{proof}

\begin{lemma}
Suppose that $g\in \Z[x)$ and that $g\in\text{ker}(e_P)$ for all stable
non-archimedean $P$. Then $g$ is divisible by $f$ in $\Z[x)$.
\end{lemma}
\begin{proof}
Let $O\subset \Q(\beta)$ be the ring of algebraic integers. It is
equal to the intersection of all non-archimedean $O_P$. The two
previous lemmas imply that $g=(\beta-x)h$ with $h\in O[x)$. Let
$K_f$ be the splitting field of $f$ over $\Q$ and let $\beta_i$ be
a conjugate of $\beta$ in $K_f$. By applying
Lemma~\ref{stableinverse} to $\Q(\beta_i)$ we see that
$g=(\beta_i-x)\sigma(h)$ in $O_f[x)$, the ring of power series
over the algebraic integers in $K_f$. So  in $O_f(x)$, $g$ is
divisible by all conjugates $x-\beta_i$. Since
$f=c(x-\beta_1)\ldots(x-\beta_d)$ for some integer $c$, this
implies that $g=(f/c) m$ for $m\in O_f[x)$. On the other hand,
since $g$ and $f$ are integral, $g/f\in \Q[x)$, and so
$m\in\Z[x)$. Since $f m=0\text{ mod } c$ and since $f$ is
primitive, $m=0\text{ mod }c$. Hence $m/c\in\Z[x)$.
\end{proof}

\begin{lemma}
Suppose that the associated polynomial $f$ is irreducible.
Then \[\zrightsup\supset (f)=\bigcap\{\text{ker}(e_P)\colon P\in \St\}.\]
\end{lemma}
\begin{proof}
It is obvious that $(f)$ is contained in the intersection of these
kernels, so we only need to prove that the opposite inclusion
holds. Suppose that $e_P(g)=0$ for all $P\in\frak S$. By the
corollary above $\frac g f\in\Z[x)$. As a complex function, $g$ is
holomorphic on the unit disc and it has zeroes at the stable
roots. Therefore $\frac g f\in\Z_{ac}[x)$ is holomorphic on the
unit disc.
\end{proof}

We shall say that the product of the evaluation maps $e_P$ for $P\in\St$
is the \textit{stable evaluation}, since it evaluates an element of $\zrightsup$
at all the stable roots of the associated polynomial.
We denote the stable evaluation by $e$.

\begin{theorem}
Suppose that the associated polynomial is irreducible. Then the
stable evaluation is a topological isomorphism.
\end{theorem}
\begin{proof}
We identify $\Q(\beta)$ with its image in $\prod\{K_P\colon
P\in\St\}$. Let ${\frak S}_\infty\subset\St$ be the subset of
archimedean primes and let ${\frak S}_p\subset\St$ be the subset
of non-archimedean primes. Let $K_\infty=\prod_{P\in{\frak
S}_\infty} K_P$ and $K_p=\prod_{P\in{\frak S}_p} K_P$. Both
$K_\infty$ and $K_p$ are locally compact. \medbreak We argue that
the image of $e$ projects onto $K_\infty$. The image of $e$
contains $\Z[\beta,\beta^{-1}]$. The algebraic integers in
$\Z[\beta,\beta^{-1}]$ form a lattice in the product of $K_P$ if
$P$ ranges over the archimedean prime-divisors. So the image of
$\polynom$ under $e$ contains a lattice $L$ of full rank in
$K_\infty$. Let $U$ be a fundamental domain for $L$ and let $F$ be
a finite set of lattice points such that the translates $F+U$
cover $\beta^{-1}\cdot U$.

For an arbitrary $x\in K_\infty$ we have to show that there exists
a $g\in \zrightsup$ such that $e(g)$ projects onto $x$. By the
definition of $U$ there is a $v\in L$ such that $u_0=x-v\in U$. If
$u_0=0$ then $x\in L$ and since $L$ is in the projection of the
image of $e$, we are done. If $u_0\not=0$ then continue dividing
$u_0$ by $\gamma$ until $u_0/\beta^{n_0}\not\in U$, which is
possible since $\beta$ is stable. Then $u_0/\beta^{n_0}\in U+f_0$
for some $f_0\in F$. Define $u_1=u_0/\beta^{n_0}-f_0$. Continue by
induction and put $u_{i+1}=u_i/\beta^{n_i}-f_i$ and truncate the
sequence if $u_{i+1}=0$. Then
$x=v+\beta^{n_0}f_0+\beta^{n_1}f_1+\ldots$ is a sum, possibly
finite, that converges to $u_0$. Since $F\subset L$ each $f_i$ is
a projection onto $K_\infty$ of $e(p_i)$ for some polynomial
$p_i$. Since $F$ is finite there exists a $g\in\zrightsup$ such
that $e(g)$ projects onto $x$.

\medbreak
We argue that the image of $e$ projects onto $K_p$.
Let $R_0\subset\Z[\beta,\beta^{-1}]$ be the subring
of algebraic integers and let $\gamma=k\beta$ be
a multiple of $\beta$ such that $\gamma\in R_0$.
Let $R_k=\gamma^k R_0$. Then
$R_0\supset R_1\supset\ldots$
forms a descending chain that intersects in $\{0\}$
since $\gamma^k$ converges to zero.
Let $F\in R_0$ be a finite set of representatives of $R_0/R_1$.
The closure of ${R_0}$ in $K_p$ is closed and open.
For an arbitrary $x\in K_p$ we have to show that there exists
a $g\in \zrightsup$ such that $e(g)$ projects onto $x$.
Multiply $x$ by a power of $\gamma$ such that $u_0=\gamma^k x$ is
in the closure of $R_0$. There exists an $f_0\in F$ such that
$u_0 = \gamma u_1 +f_0$ and we can construct an infinite sum
that converges to $u_0$ in the same way as in the archimedean case.
\medbreak

So the image of the stable evaluation $e\colon\zrightsup/(f)\to
K_\infty\times K_p$ projects onto both factors. By local
compactness and by Baire's property, both projections are open,
and thus $e$ is an open map. Therefore, the image $I$ of the
stable evaluation is a locally compact subgroup of $K_\infty\times
K_p$, hence $I$ is closed. The factor group $K_\infty\times K_p/I$
is isomorphic to $K_\infty/(I\cap K_\infty)$, hence it is
connected. It is isomorphic to $K_p/(I\cap K_p)$ as well, hence it
is totally disconnected. So $e$ has to be a surjection. The
results above imply that the kernel of the stable evaluation $e$
is equal to $(f)$.
\end{proof}

Summarizing these results,
the stable group $\Omega^+$ is isomorphic to $\prod\{K_P\colon P\in\St\}$ and
the isomorphism is induced by the evaluation of a power series at the stable roots
of the associated polynomial. By symmetry $\Omega^-$ is isomorphic to
$\prod\{K_P\colon P\in\Un\}$ under the evaluation at the unstable roots.
So algebraically, the exact sequence \ref{unfolding} is an embedding of
the ring $\Z[\beta,\beta^{-1}]$ in a product of completions of $\Q(\beta)$,
such that its cokernel is its Pontryagin dual.
In other words, it is the standard algebraic discretization of a ring.
In the dynamic setting, it is not easy to find a fundamental domain for
$\Omega$ in $\Omega^-\times\Omega^+$, but in the algebraic setting,
it is.

\begin{theorem}
Let $R$ be the ring of algebraic
integers in $\Z[\beta,\beta^{-1}]$.
Let $K_a$ be the product of the archimedean completions of $\Q(\beta)$.
Let $K_{na}$ be the product of the non-archimedean completions that are
stable or unstable.
Let $U$ be a fundamental domain of $R$ in $K_a$ and let $O_{na}$
be the ring of integers in $K_{na}$. Then
$O_{na}\times U$ is a fundamental domain of $\Z[\beta,\beta^{-1}]$
in the unfolding.
\end{theorem}
\begin{proof}
Assume for the moment that $\Z[\beta,\beta^{-1}]$ is dense in
$K_{na}$. Then for every $\kappa\in K_a\times K_{na}$
there exists a $q\in\Z[\beta,\beta^{-1}]$ such that $\kappa-q\in
K_a\times O_{na}$. By translation of an $r\in R$ we can move $\kappa-q$
to an element $v\in U\times O_{na}$. The element $v$ is unique,
for if $v-v'\in \Z[\beta,\beta^{-1}]$ for some $v'$ then $v-v'\in
R$ since $v-v'\in O_P$ for all non-archimedean $P$ that are stable
or unstable and since $\Z[\beta,\beta^{-1}]\subset O_P$ for all
other non-archimedean prime-divisors.

It remains to show the validity of our assumption. By the
approximation theorem it suffices to show that
$\Z[\beta,\beta^{-1}]$ is dense in $\Q(\beta)\subset K_{na}$. The
closure of the ring of algebraic integers in $\Z[\beta,\beta{-1}]$
is open in $\Q(\beta)$. So the closure of $\Z[\beta,\beta^{-1}]$
contains a neighborhood of $0$. The sequence
$1/(\beta^{n}+\beta^{-n})$ converges to $0$. Hence, for every
$x\in \Q(\beta)$ there exists an $n$ such that
$x/(\beta^{n}+\beta^{-n})$ is in the closure of
$\Z[\beta,\beta^{-1}]$.
\end{proof}

So far we only considered irreducible associated
polynomials, but it is not difficult to extend the result to
general associated polynomials. By the Chinese remainder theorem
the stable group is isomorphic to the direct sum of
$\zrightsup/(f_i)$ over all primary factors $f_i$ of the
associated polynomial, so we need to extend the results to primary
polynomials only.

\begin{theorem}
Suppose that the associated polynomial is primary $f=g^n$ for an
irreducible polynomial $g$. Then $\zrightsup/(f)$ is isomorphic to
$\prod\{K_P^n\colon P\in\St\} $ as a topological group.
\end{theorem}

\begin{proof}
Let $m^{(j)}$ denote the $j$-th derivative of the power series
$m$. Obviously, if $m\in (g^n)$ then $m\in (g^{n-1})$
and $m^{(n)}\in (g)$. Conversely, if
$m=g^n\cdot m_1$ for some power series $m_1$
and if $g$ divides $m^{(n)}$, then $g$ divides $(g')^n\cdot m_1$. By
the irreducibility of $g$ there exist polynomials $h_1,h_2$ such
that $h_1(g')^n+h_2g=1$. Therefore $g$ divides $m_1$.
So $m\in (g^n)$.
By induction, $m\in(g^n)$ if and only if $m^{(j)}\in (g)$ for
all $j<n$. In particular, the map $m\mapsto
(m,m',\ldots,m^{(n-1)})$ induces an injective group homomorphism
$\varphi\colon\zrightsup/(\bar f)\to \bigoplus_{i=1}^n \zrightsup^+/(\bar g)$.
This map is also surjective. To
prove this, let $h_k$ be a polynomial such that
$h_k(g')^k=1\text{ mod } (g)$ and such that $m=h_k g^k$. Then
$m^{(i)}\in (g)$ for $i<k$ and $m^{(k)}=1 \text{ mod }(g)$, so for
a polynomial $h$, $(h\cdot m)^{(i)}\in (g)$ for $i<k$ and $(h\cdot
m)^{(k)}= h\text{ mod } (g)$.
The composition of $\varphi$ and the product of the evaluation map $e$ is clearly
continuous.
\end{proof}

\begin{corollary}\label{stable}
The associated polynomial is monic (up to a sign) if and only if
all unstable prime-divisors are archimedean.
\end{corollary}
\begin{proof}
Let $f=a_dx^d+\ldots+a_0$. By the Newton
polygon non-archimedean unstable primes-divisors extend the primes
that divide $a_d$ while the unstable prime divisors extend the primes
that divide $a_0$. So
$a_0$ is a unit if and only if $\Q(\beta)$ has no unstable
non-archimedean prime-divisors.
\end{proof}

If the expansive automorphism is defined on the torus $\T^n$, then
$a_d$ and $a_0$ are units, so all stable and unstable prime-divisors are
archimedean. The exact sequence of $\T^n$ is simply the universal
factorization through the universal covering space $\R^n$.
This is the case that has been considered
by Vershik.

\section{The path-component of the identity}

We show that if the associated polynomial $f$ is irreducible, then
the homoclinic group $\Omega$ is irreducible in the following way:
its intersection with the path-component of the identity is either
$\{0\}$ or $\Omega$.
\medbreak
We denote the union $\St\cup\Un$ by $\frak P$.
For any subset $\frak Q\subset\frak
P$ we denote the factor ring $\prod_{P\in\frak
Q}K_P\times\prod_{P\in{\frak P}-{\frak Q}}\{0\}$ by $\frak {U_Q}$.
In particular, $\Omega^+\times\Omega^-$ is identical to $\frak
{U_P}$. Since $\frak{U_Q}$ is a factor of $\frak {U_P}$ it
projects onto $\Gamma$.

\begin{theorem}\label{pathcomp}
Suppose that $(\Gamma,\alpha)$ is a cyclic expansive system such
that the associated polynomial $f$ is irreducible. Let
$\Gamma_0\subset\Gamma$ be the path-component of the identity. The
following statements are equivalent:
\begin{enumerate}
\item $\Omega\cap \Gamma_0\not=\{0\}$ \item
$\Omega\subset\Gamma_0$ \item $f(x)$ or $f(1/x)$ is monic (up to a
sign).
\end{enumerate}
\end{theorem}
\begin{proof}

$(2)\implies(1)$ is trivial.

$(1)\implies (3)$
Let $\frak Q\subset \frak P$ be the subset of archimedean prime-divisors.
Then $\frak{U_Q}$ is the path-component of
the identity of $\Omega^-\times\Omega^+$ and the projection of $\frak {U_Q}$
onto $\Gamma$ is equal to $\Gamma_0$.
Suppose that $\Omega\cap\Gamma_0\not=\{0\}$. By
the exact sequence~\ref{unfolding}, the preimage of $\Omega\subset\Gamma$
in $\Omega^-\times\Omega^+$ is equal to $\Omega\times\Omega$.
So
$\frak{U_Q}$
intersects $\Omega\times\Omega$ in a point outside the origin. This
point is represented by $(g,h)$ for two Laurent polynomials
$g$ and $h$. Let $P$ be a prime-divisor in ${\frak
P}\setminus{\frak Q}$ and to fix our ideas, suppose that $P$ is
stable. Since $(g,h)\in \frak{U_Q}$ necessarily
$h(\beta)=0\in K_P$, hence $f$ divides $h$, so $(g,h)$ is equivalent to
$(g,0)$. Then $g\not=0$ and all unstable prime-divisors are archimedean.
By Corollary~\ref{stable}, $f$ is monic up
to a sign. By symmetry, if $\frak Q$ contains all unstable
prime-divisors, then $f(1/x)$ is monic.

$(3)\implies(2)$ Suppose that $f$ is monic. Then all unstable
prime-divisors are archimedean, and so $\Omega^-\subset\Gamma_0$.
Hence $\Omega\subset\Gamma_0$.
\end{proof}

In the proof of the theorem we have used only that $\frak Q$ is a
set of archimedean prime-divisors. So the proof remains valid if
$\frak Q$ is a subset of the archimedean prime-divisors. In other
words, if the image of $\frak U_Q$ contains a non-trivial
homoclinic point, then $\frak Q$ contains all stable or all
unstable prime-divisors. This corollary of the proof is a form of
the Pisot-Vijayaraghavan theorem.

\begin{corollary}[weak form of the Pisot-Vijayaraghavan Theorem, cmp.~\cite{Cassels}]
Suppose that $\beta>1$ is a real algebraic number and suppose that
$\beta$ is hyperbolic. If for some $t\in \R$ the sequence
$t\beta^n\text{ mod }1$ converges to $0$ as $|n|\to \infty$. Then
$\beta$ is an algebraic integer and the absolute value of all
conjugates of $\beta$ is $<1$.
\end{corollary}
\begin{proof}
Let $f$ be the minimum polynomial of $1/\beta$. The formal power
series $g=\sum (t\beta^n)x^n\in\T(x)$ is annihilated under
multiplication by $f$. So we may consider $g$ as a non-trivial
homoclinic point in the compact abelian group $\Gamma$ with
associated polynomial $f$. Even more so, $g$ is contained in the
arcwise-connected subgroup $\{(r\beta^n)_{n\in \Z}\colon r\in\R\}$
of $\Omega^-$. By the theorem above $\beta$ is an algebraic
integer.

Let $\lfloor x\rfloor $ denote the integer part of the real number
$x$. A preimage of $g$ in $\R_{ac}(x)$ is given by
\[\sum \left((t\beta)^n-\lfloor (t\beta)^n\rfloor\right)x^n.\]
Under multiplication by $f$ it is mapped onto $-f\cdot\sum \lfloor
(t\beta)^n\rfloor x^n$. Note that $\lfloor t\beta^n\rfloor=0$ if
$n$ is sufficiently small. One verifies that the coefficients of
the formal power series
\[(x-1/\beta)\cdot \sum \lfloor t\beta^n\rfloor x^n=\sum
\frac{\beta\lfloor t\beta^n\rfloor-\lfloor
t\beta^{n+1}\rfloor}{\beta}x^{n+1}\] are bounded by
$\max\{1,|t|\}$, so it represents a holomorphic function in the
unit disc. In particular, $-f\cdot\sum \lfloor t\beta^n\rfloor
x^n$ is equal to zero at all stable roots of $f(1/x)$ other then
$1/\beta$. In particular, there is only one stable prime-divisor
$P$ of $\Q(1/\beta)$ for which $e_P(-f\cdot\sum \lfloor
t\beta^n\rfloor x^n)\not=0$. Hence, if we use $\{P\}$ for $\frak
Q$, then we see that the image of $\frak{U_Q}$ contains a
non-trivial homoclinic point. So, as a corollary of the proof of
the previous theorem, $P$ is the only stable prime-divisor.
\end{proof}

Our corollary is a weak form of the Pisot-Vijaraghavan Theorem since
we need the additional (and superfluous) assumption that $\beta$
is hyperbolic. The Pisot-Vijayaraghavan Theorem also gives an explicit
description of $t$, but we could have given an explicit
description as well, since $(t\beta^n)$ is in the homoclinic group.
\medbreak
An algebraic integer $\beta$ for which all conjugates have absolute
value $<1$ is called a Pisot-Vijayaraghavan number, or simply a
Pisot number.

\section{Symbolic codings}

By  a \textit{coding} of an automorphism $\alpha$ of $\Gamma$, we mean a
continuous mapping from a
closed, shift invariant subset $S$  of $\{ 0, \dots, N  \}^\Z $
onto $\Gamma$ commuting
the shift map and $\alpha$. We call an element of $S $  \textit{finite} if all but
finitely many terms are $0$. Arithmetic codings, which have the property that homoclinic
points are coded by finite symbolic sequences, seem to have
first appeared in \cite{Bert} and
were more fully developed by Vershik \cite{Vershik} and Sidorov \cite{Sidorovvershik}
\cite{Sidorov1}.
One may verify that if a coding is arithmetic, then the stable group and the unstable
group are coded by symbolic sequences that have a tail of zeroes.

\medbreak Let $\beta>1$ be a real number. The one-sided
$\beta$-shift $Z_\beta \subset \{ 0, \dots, \lfloor \beta \rfloor
\}^\N $ is conjugate to $x\mapsto \beta x\text{ mod }1$ under the
projection $(s_n)\mapsto \sum_{n=1}^\infty s_n\beta^{-n}$ , see
\cite{Blanchard} for a survey.   If $\beta$ is a Pisot number,
then $Z_\beta$ is a sofic set and the projection $Z_\beta\to[0,1]$
is $1-1$ on the complement of a countable set.  The two-sided
$\beta$-shift $S_\beta$ is obtained by augmenting every element of
$Z_\beta$ with a left tail of zeroes and taking the closure of the
orbits under the shift. If $\beta$ is a Pisot unit, its minimal
polynomial is associated to an expansive toral automorphism, that
is naturally coded by $S_\beta$. The coding is given by
\[
(s_n)_{n\in\Z}\mapsto \lim_{N \to \infty} \left(\sum_{n=-N}^{n=N}
 s_n \beta^{-n} \right)\bold{ e} \text{ mod }  \Z^n
\]
where $\bold{ e}$ is a point of $\R^n$ representing a homoclinic
point of the automorphism. Vershik and Sidorov gave a detailed
analysis of this coding, analyzing the effect of the choice of the
homoclinic point $\bold{ e}$, which amounts to a choice of
lattice.

\medbreak A symbolic coding $S\to\Gamma$ is called \textit{almost}
$1-1$ if there exists a subset $\Gamma_0\subset\Gamma$ of full
Haar measure such that the fibers over $\Gamma_0$ are singletons.
If the polynomial associated to an automorphism is irreducible
with a Pisot root  $\beta$, Schmidt~\cite{Schmidt} showed that
$S_\beta$ codes the automorphism bounded-to-one,  and  for a class
of Pisot units he was also able to show that the coding is almost
$1-1$. This led to his conjecture that one could always find an
almost $1-1$ coding in the case of a Pisot unit. Sidorov and
Vershik found that for a wider class of Pisot units (including all
quadratic units), $S_\beta$ gives an almost $1-1$ coding of toral
automorphisms. Sidorov has reduced the conjecture to algebraic
considerations ~\cite{Sidorov1},\cite{Sidorov2}.

\medbreak We call  an automorphism $ \alpha \sim f(x)$  a
\textit{Pisot automorphism} if $f(x)$ or $f(1/x)$ is (up to a
sign) monic and irreducible with a Pisot root $\beta$. And in this
context, $\beta$ is the \textit{associated Pisot number}. If an
automorphism is not Pisot, then one cannot expect $S_\beta$ to
code the automorphism almost $1-1$ since the entropies in these
cases do not necessarily match. For example, consider $f(x)=2x-3.$
The entropy of $S_{3/2}$ is $\text{log} ( 3/2 )$,  while the
associated automorphism has entropy $ \text{log} (3)$ by
Yuzvinskii's formula, (see, e.g., \cite{LW}).

One is naturally led to the following generalization of Schmidt's
conjecture.

\begin{conjecture}
A  Pisot automorphism with associated Pisot number $\beta$ admits
an almost $1-1$ coding by $S_\beta$.
\end{conjecture}

Recall that we are only considering cyclic automorphisms: the
above conjecture is false in general even for tori. Although we
have not been able to prove this conjecture, we can establish it
in certain cases that generalize previous results for Pisot units.
As in \cite{FS}, let $\text{Fin}( \beta )$ denote the set of those
numbers $x \geq 0 $ admitting a finite $\beta$ expansion. We call
$\beta$ \textit{finitary} if $\text{Fin} ( \beta )
=(\Z[\beta^{-1}])_+.$   The Pisot root of a polynomial of the form
$f(x)=x^n-a_{n-1}x^{n-1}+ \cdots - a_0$ with $a_{n-1} \geq a_{n-2}
\geq \cdots \geq a_0 \geq 1$ is finitary \cite{FS}, but it is not
clear which other $\beta$ are finitary.
\begin{theorem}
If $\alpha$ is a  Pisot automorphism with an associated finitary
Pisot number $\beta$, then the natural coding $S_\beta \to \Gamma$
\[
(s_n)_{n\in\Z} \mapsto \sum_{-\infty}^{\infty} s_n x^n +(f)
\]
is almost $1-1$.
\end{theorem}

\begin{proof}
We treat the case that $f(1/x)$ is monic up to a sign, the other
case being handled in a similar way. In this case, the stable
manifold $\Omega^+$ has just one stable divisor, and it is
archimedean. We identify the stable manifold with  $\R$ via the
stable evaluation $e$, which then conjugates the action of the
automorphism on $\Omega^+$ with multiplication by $\beta$. If we
restrict $e$ to those elements of $S_\beta$ with a left tail of
zeros, we obtain the ``positive'' ray in the stable manifold. We
now use the unfolding of Theorem~\ref{unfolding} to analyze the
coding.  For $ (s_n) \in S_\beta$, let $c^+ ( (s_n)
)=\sum_{n=1}^{\infty} s_n x^n +(f) \in \Omega^+$ and $c^- ( (s_n)
)=\sum_{- \infty}^{n=0} s_n x^n +(f) \in \Omega^-$. Then the
natural coding corresponds to adding $c^+ ( (s_n) )$  and $c^- (
(s_n) )$ in $\Gamma$.

Let $S_{\beta}^{-}$ consist of those elements $(s_n) \in S_\beta$
for which $s_n=0$ for $n>0$. For a class of Pisot numbers
including the finitary, it has been shown (see, \cite{Sol},
\cite{Grab}, \cite{FS}) that $S_{\beta}^{-}$ admits a group action
that extends (measure-theoretically) the odometer or adic
transformation. This transformation maps points with an infinite
tail of zeros to the left to their immediate successor in the
lexigraphic ordering. Notice that when $\beta$ is an integer, this
odometer map corresponds to addition by $1$ when one identifies
sequences with their evaluations. In general it does not
correspond to addition by $1$ in this sense: sequences in
$S_{\beta}^{-}$ with a left tail of zeroes do not typically
evaluate to an integer. However, when $\beta$ is finitary this map
extends to $S_{\beta}^{-}$ as a minimal map with purely discrete
spectrum. We now relate this map to the return map of a flow on
the two sided shift $S_\beta$ which has purely discrete spectrum
and whose associated action extends the addition of evaluations of
sequences with left tails of zeros.

\medskip For finitary $\beta$, the ``carrying over'' of digits when adding
two elements of $\text{Fin}( \beta )$  has a uniform bound~
\cite[Proposition 2]{FS}. Thus, there is a uniform gap length $G$
so that two elements of $S_\beta$ which have infinitely many (to
the right and to the left) matching blocks of $0$'s of length at
least $G$ can be added by treating each pair of sequences between
successive blocks of $0$'s as an element of $\text{Fin}( \beta )$.
The sum then  corresponds to the sequence formed by concatenating
all the resulting finite expansions of sums, see
\cite{Sidorovvershik}. This group structure then coincides with
addition of expansions and is defined on a subset of $S_\beta
\times S_\beta$  of full measure, provided only that the measure
is shift-invariant and positive on cylinders. One can define an
additive action of $\R$ (a flow)  on $S_\beta$ by identifying
points in $\R$ with their greedy $\beta$-expansion and adding
according to this group structure. We shall call this the
\textit{additive flow}, which is defined on a set of full measure
in $\R \times S_\beta$. The orbits of the additive flow are dense
and the group operation is compatible with the flow. Hence, the
additive flow has purely discrete spectrum. Under the natural
coding, the additive flow commutes with the continuous action in
$\Gamma$ given by $(t,g+(f)) \mapsto e^{-1}(t)+g+(f)$, which has
the stable manifold and its translates as orbits. (Here
$e^{-1}(t)$ is understood to mean the point in the stable manifold
that evaluates as  $t$.) Since the group structure is compatible
with this flow, it also has purely discrete spectrum. Notice that
the image under the unstable evaluation of $c^- (S_\beta)$  then
is a compact cross-section of this flow and the return map to this
cross-section coincides with the odometer or adic transformation.

\medskip Thus, the natural coding is a (measure theoretic) homomorphism of
one compact group onto another. Hence, the coding must be
$k-\text{to}- 1$ for some finite $k$ off a set of measure 0. The
collection of points of $S_\beta$ corresponding to the kernel of
the coding are shift invariant and finite in number. Thus, each
such point can be represented by a periodic sequence. Let $(p_n)$
then be any periodic sequence in the kernel. Recall that all
points in the kernel of
\[
 \Omega^-\times\Omega^+ \to \Gamma
\]
from the unfolding have homoclinic points in each factor. In terms
of the coding, this means that the stable evaluation of $c^+(
(p_n))$ must then be a point in $ ( \Z[\beta^-])_+ = \text{Fin}(
\beta ) $. Thus, the positive portion of  $(p_n)$ must have the
same evaluation as a number with a finite $\beta$-expansion. The
only alternative admissible expansion of a finite sequence is a
sequence ending with the right tail of the infinite expansion of
1. Thus, if $(p_n)$ is not the zero sequence, it must be a
sequence each right tail of which coincides with a tail of the
infinite expansion of $1$. But each such sequence (by
considerations of continuity) must represent the identity of the
group.
\end{proof}

It might be objected that the above described coding is abstract.
However, by identifying the polynomial $1$ with a specific
fundamental homoclinic point of $\Gamma$, one can easily form a
geometric interpretation of the coding. The additive flow maps the
identity to the chosen homoclinic point in one unit of time, and
the remainder of the coding is determined by continuity.

\begin{example}
If $f(x)=2x-1$, then $\beta=2$ and $S_\beta= \{0,1 \}^\Z$. In this
case 2 is the only unstable prime divisor, the additive flow is
the suspension of the standard adding machine on the $2-adic$
integers, which form the unstable cross-section of the flow. The
zero element of  $S_\beta$ is also represented by the sequence of
all $1$'s.
\end{example}

\begin{example}
If $f(x)=x^2+nx-1$ $(n>1);$ $f(1/x) \approx x^2-nx-1 $. The
additive flow is an irrational flow on a torus and the return time
to the unstable cross-section is a step function with two steps.
The zero element of $S_\beta$ is also represented by the periodic
sequences composed of $n0$'s.
\end{example}

\begin{example}
If $f(x)=2x^2+nx-1$ $ ( n>2)$;  $ f(1/x) \approx x^2-nx-2 $. Then
$S_\beta$ is a subshift of finite type and there is both an
archimedean and a non-archimedean unstable divisor. The unstable
cross-section is one-dimensional but disconnected, and the return
map  to the cross-section is the product of an interval exchange
and a $2-adic$ adding machine. The zero element of $S_\beta$ is
also represented by the periodic sequences composed of $n1$'s.
\end{example}

\section{Acknowledgement}

This paper would not have been written without the help of Hendrik
Lenstra. \newline This paper was written while the first author
was visiting Delft University, supported by an NWO visitor's
grant.

\end{document}